\documentclass[a4paper,12pt]{article}

\usepackage{amsmath,amssymb,amsthm,amsfonts,epsfig}

\numberwithin{equation}{section}

\newtheorem{thm}{Theorem}[section]
\newtheorem{lem}{Lemma}[section]
\newtheorem{rem}{Remark}[section]
\newtheorem{prop}{Proposition}[section]
\newtheorem{cor}{Corollary}[section]

\title{Asymptotic Behavior of Blowup Solutions for Elliptic Equations with
Exponential Nonlinearity and Singular Data}

\author{
 {\sc Lei Zhang} \\
 {\sc \footnotesize  University of Alabama at Birmingham}  \\
 {\sc \footnotesize
    Department of Mathematics,}\\
    {\sc \footnotesize 452 Campbell Hall} \\
 {\sc \footnotesize
    1300 University Boulevard}\\
 {\sc \footnotesize
    Birmingham, AL 35294-1170}  \\
{\tt \footnotesize
    leizhang@math.uab.edu}\thanks{Supported by
National Science Foundation Grant 0600275. Running title:
Asymptotic expansion for blowup solutions} }

\date{}

\input { amssym.def}

\begin{document}
\maketitle

\begin{abstract}
We consider a sequence of blowup solutions of a two dimensional,
second order elliptic equation with exponential nonlinearity and
singular data. This equation has a rich background in physics and
geometry. In a work of Bartolucci-Chen-Lin-Tarantello it is proved
that the profile of the solutions differs from global solutions of
a Liouville type equation only by a uniformly bounded term. The
present paper improves their result and establishes an expansion
of the solutions near the blowup points with a sharp error
estimate.
\end{abstract}

\noindent

{\bf Mathematics Subject Classification (2007):
 35J60, 35B45, 53C21}

 \smallskip

 {\bf Keywords:}
 Liouville equation, Blowup analysis.

\section{Introduction}

Two dimensional semilinear elliptic equations with exponential
nonlinearities arise naturally in conformal geometry and physics.
The study of these equations is always related to their blowup
phenomena. When a sequence of solutions tends to infinity near a
blowup point, the asymptotic behavior of the solutions near the
blowup point carries important information. In some applications
it is crucial to completely understand the asymptotic behavior of
blowup solutions. In this article we study the following equation:
\begin{equation}
\label{jan2e1} \Delta u+|x|^{2\alpha}H(x)e^u=0,
 \quad \hbox{in } B_1\subset \mathbb R^2,
\end{equation}
where $B_1$ is the unit ball in $\mathbb R^2$, $\alpha \in \mathbb
R^+ \setminus \mathbb N,$ ($\mathbb N$ is the set of natural
numbers) and $H \in C^3 ({\overline B_1})$ is a positive function.
If a sequence of solutions $\{u_i\}$ tends to infinity near a
point other than $0$, the coefficient function $|x|^{2\alpha}H(x)$
is bounded above and below near the blowup point. This situation
has been extensively studied and the asymptotic behavior of the
blowup solutions is well understood ( see, for example \cite{BT},
\cite{BM},\cite{Chenxx},\cite{Licmp},\cite{LiSha},\cite{ChenLin1},\cite{zhangcmp}).
In this article we mainly consider a sequence of solutions
$\{u_i\}$ of (\ref{jan2e1}) such that
\begin{equation}\label{mar28e3}
  u_i(z_i)= \max_{B_1} u_i \to \infty,
  \quad z_i\to 0,
  \quad 0 \mbox{ is the only blowup point in } {\overline B_1}.
\end{equation}
 We shall describe the asymptotic profile of $\{u_i\}$ near $0$ under natural
assumptions on $H$ and the oscillation of $\{u_i\}$ on $\partial
B_1$.

The blowup analysis for (\ref{jan2e1}) near $0$ reflects the
bubbling feature of a few important equations or systems of
equations in physics or geometry. For example, the following mean
field equation is defined on Riemann surfaces:
\begin{equation}
\label{mar29e1} \Delta_gw+\rho \left(
\frac{h(x)e^w}{\int_Mh(x)e^wdV_g}-\frac 1{|M|} \right)
=4\pi\sum_{j=1}^m\alpha_j(\delta_{p_j}-\frac 1{|M|})
\end{equation}
where $M$ is a compact smooth Riemann surface without boundary,
$h$ is a positive smooth function on $M$, $\rho$ is a positive
constant, $|M|$ is the volume of $M$, $\Delta_g$ is the
Laplace-Beltrami operator and  $\alpha_j\delta_{p_j}$ are Dirac
sources. For (\ref{mar29e1}), the profile of blowup solutions near
each $p_j$ is exactly described by (\ref{jan2e1}) and $\alpha_j$
in (\ref{mar29e1}) plays the same role as $\alpha$ in
(\ref{jan2e1}). From the physical point of view, it is important
to consider the case $\alpha_j> 0$ in (\ref{mar29e1}) as it is
closely related to the self-dual equations in the Abelian
Chern-Simons-Higgs theory (see
\cite{Dunne},\cite{Hong},\cite{Jackiw}, etc) and the Electroweak
theory (see \cite{ambjorn}\cite{lai}, etc). There is a
considerable literature on the mean field equation (\ref{mar29e1})
and closely related topics, we list the following as a partial
list:
\cite{BT}\cite{Chang}\cite{ChenLin1}\cite{ChenLin2}\cite{ChenLin3}\cite{ding}
\cite{esposito}\cite{lincs1}\cite{lincs2}\cite{linlucia1}\cite{lucia3}\cite{lucia1}\cite{lucia2}.

Another application comes from the 2-dimensional open Toda system
for $SU(N+1)$. The simplest example is
$$\left\{\begin{array}{ll}
-\Delta u_1^k=2h_1^ke^{u_1^k}-h_2^ke^{u_2^k}\\
\\
-\Delta u_2^k=2h_2^ke^{u_2^k}-h_1^ke^{u_1^k}
\end{array}
\right. \quad \mbox{in}\quad B_1
$$
where $u_1^k$, $u_2^k$ are sequences of blowup solutions and
 $h_1^k$, $h_2^k$ are positive, smooth functions very close to $1$.
Even for this simplest example, the blowup analysis is challenging
because $u_1^k$ and $u_2^k$ may have common blowup points and the
asymptotic behavior of them near their common blowup points is not
yet well understood (see \cite{jostwang},\cite{jostwang2},
\cite{linwang},\cite{japan},\cite{lucia1} and the references
therein for recent development). The analysis for (\ref{jan2e1})
is closely related to the Toda system and the result we prove in
this work (Theorem \ref{thm2}) helps to understand the system.

In addition to the background in physics, (\ref{jan2e1}) has a
well known interpretation in geometry. Let $g_0$ be the Euclidean
metric on $B_1$, then $\frac 12|x|^{2\alpha}H(x)$ is the Gauss
curvature under metric $e^u g_0$. In this sense (\ref{jan2e1}) is
related to the Nirenberg problem or more generally the
Kazdan-Warner problem.

When $\alpha$ is equal to $0$, the behavior of a sequence of
blowup solutions $\{u_i\}$ to
\begin{equation}
\label{mar28e1} \Delta u_i + H_i(x) e^{u_i}
 =0
 \quad \mbox{in}\quad B_1
\end{equation}
has been extensively studied through the works of Brezis-Merle
\cite{BM},
 Li-Shafrir\cite{LiSha},
Li\cite{Licmp}, Chen-Lin\cite{ChenLin1} and the references
therein. In \cite{Licmp} Li proved that if a sequence of blowup
solutions $\{u_i\}$ of (\ref{mar28e1}) has a bounded oscillation
near their blowup point, then in a
 neighborhood of this point $\{u_i\}$ is only $O(1)$ different from a sequence of
standard bubbles appropriately scaled. Later Chen-Lin
\cite{ChenLin1} and the author \cite{zhangcmp} improved Li's
estimate to the sharp form by different approaches. Since the case
$\alpha>0$ is more meaningful in physics, it is important to
obtain similar results for equation (\ref{jan2e1}) when
$\alpha>0$.

In this article we address the case $\alpha\not \in \mathbb N$.
Our work is based on a result by Bartolucci-Chen-Lin-Tarantello
\cite{BCLT} who studied (\ref{mar29e1}) and established a result
similar to Li's result for $\alpha=0$. More specifically, let
$\{u_i\}$ be a sequence of functions solving
\begin{equation}
\label{mar28e2} \Delta u_i+|x|^{2\alpha}H_i(x)e^{u_i}=0,
 \quad \hbox{in } B_1,
\end{equation}
 such that \eqref{mar28e3} holds.
 Suppose $\{u_i\}$ has bounded oscillation on $\partial B_1$ and a
 bound on the energy:
\begin{equation} \label{mar28e4}
  \left\{
  \begin{array}{c}
  |u_i(x)-u_i(x')|\le C_0 \quad \forall x, x'\in
\partial B_1,
  \\\vbox{\vskip6mm}
  \int_{B_1}|x|^{2\alpha}H_i(x)e^{u_i}\le C_0.
  \end{array}
  \right.
\end{equation}
Then the following result is proved in \cite{BCLT}:
\begin{thm}(Bartolucci-Chen-Lin-Tarantello)
\label{thm1} Let $\{u_i\}$ satisfy (\ref{mar28e2}),
(\ref{mar28e3}), (\ref{mar28e4}) and let $H_i$ satisfy
$$
  0<\frac 1{C_1}\le H_i(x)\le C_1, \quad \|\nabla
 H_i\|_{L^{\infty}(B_1)}\le C_1.
$$
Then for $\alpha \in \mathbb R^+ \setminus \mathbb N$, there
exists $C>0$ such that
$$
 \left|
   u_i(x)-\log
 \frac{e^{u_i(0)}}{(1+\frac{H_i(0)}{8(1+\alpha)^2}e^{u_i(0)}|x|^{2\alpha+2})^2}
 \right|\le C \quad \hbox{in } B_1.$$
\end{thm}
Theorem \ref{thm1} is a refinement of a result in
Bartolucci-Tarantello \cite{BT}, concerning the quantization
phenomena for blowup solutions of (\ref{jan2e1}) by using an
improved argument based on the Pohozaev Identity. The idea of
using Pohozaev's type arguments was first introduced by
Bartolucci-Tarantello \cite{BT} in the analysis of the bubbling
phenomena for (\ref{jan2e1}) with $\alpha>0$.

In our main result below we assume the following on $H_i$:
\begin{equation}
\label{nov6e1}
  0<\frac 1{C_1}\le H_i(x)\le C_1, \quad
  \|H_i\|_{C^3(B_1)}\le C_1,
\end{equation}
 and we consider the harmonic function $\psi_i$ solving
\begin{equation}
\label{nov27e2} \left\{\begin{array}{ll}
\Delta \psi_i=0  \quad \hbox{in } B_1,\\
\\
 \psi_i=u_i-\frac{1}{2\pi}\int_{\partial B_1}u_idS
 \quad \hbox {on } \partial B_1.
\end{array}
\right.
\end{equation}
Clearly $\psi_i$ is a bounded function on $B_1$ and $\psi_i(0)=0$.
 For simplicity, we define
$$
  V_i(x)=H_i(x)e^{\psi_i(x)},
$$
and introduce the following two constants:
\begin{equation}
\label{mar26e1}
 \Lambda_1=-\frac{\pi}{V_i(0)
 \sin \left(\frac{\pi}{1+\alpha} \right)(1+\alpha)}
 \left( \frac{8(1+\alpha)^2}{V_i(0)} \right)^{\frac{1}{1+\alpha}},
\end{equation}
\begin{equation}
\label{mar26e2}
\Lambda_2=\frac{\pi}{V_i^2(0)\sin(\frac{\pi}{1+\alpha})(1+\alpha)}
 \left( \frac{8(1+\alpha)^2}{V_i(0)} \right)^{\frac{1}{1+\alpha}}.
\end{equation}
 Note that $\psi_i(0)=0$ implies the following:
\begin{eqnarray*}
V_i(0)=H_i(0),\quad \nabla V_i(0)=\nabla H_i(0)+H_i(0)\nabla
 \psi_i(0)\\
\Delta V_i(0)=\Delta H_i(0)+2\nabla H_i(0)\cdot \nabla
 \psi_i(0)+H_i(0)|\nabla \psi_i(0)|^2.
\end{eqnarray*}
Our main theorem improves Theorem~\ref{thm1} as follows:
\begin{thm}
\label{thm2}
 Let $\alpha\in \mathbb R^+ \setminus \mathbb N$. Suppose $u_i$
 satisfies (\ref{mar28e2}), (\ref{mar28e3}), (\ref{mar28e4}) and
(\ref{nov6e1}). Then in $B_1$ there holds
\begin{eqnarray*}
u_i(x)&=&\log
 \frac{e^{u_i(0)}}{(1+\frac{V_i(0)}{8(1+\alpha)^2}e^{u_i(0)}|x|^{2\alpha+2})^2}
+\psi_i(x)\\
&&-\frac{2(1+\alpha)}{\alpha V_i(0)} \frac{\nabla V_i(0)\cdot
 x}{1+\frac{V_i(0)}{8(1+\alpha)^2}e^{u_i(0)}|x|^{2\alpha+2}}\\
&&+\bigg (\Lambda_1\Delta V_i(0) +\Lambda_2|\nabla V_i(0)|^2\bigg
)\log \left( 2+e^{\frac{u_i(0)}{2(1+\alpha)}}|x| \right)
e^{-\frac{u_i(0)}{1+\alpha}} +O(e^{-\frac{u_i(0)}{1+\alpha}}),
\end{eqnarray*}
 where $\psi_i$, $\Lambda_1$, $\Lambda_2$ are defined by
(\ref{nov27e2}), (\ref{mar26e1}) and (\ref{mar26e2}),
respectively.
\end{thm}
Note that we use $O(e^{-\frac{u_i(0)}{1+\alpha}})$ to denote a
smooth
 function in $B_1$ whose
$C^3$ norm is bounded by
 $C(\alpha,C_0,C_1)e^{-\frac{u_i(0)}{1+\alpha}}$.

Theorem \ref{thm2} corresponds to the results of Chen-Lin
\cite{ChenLin1} and the author \cite{zhangcmp} for the case
$\alpha=0$. One important application of getting sharp blowup
estimates is to derive a degree counting formula for the mean
field equation (\ref{mar29e1}). When the right hand side of
(\ref{mar29e1}) is $0$, Chen-Lin \cite{ChenLin2} established a
degree counting formula in terms of the genus of the Riemann
surface using the sharp estimate in \cite{ChenLin1}. Another way
of counting the degree has been proved by Malchiodi
\cite{Malchiodi}. We expect to use Theorem \ref{thm2} to derive a
degree counting formula for the general mean field equation
(\ref{mar29e1}) in a forthcoming paper.

The assumption $\alpha\not\in \Bbb N$ is essential to Theorem
\ref{thm2}. When (\ref{mar28e2}) is compared with the case
$\alpha=0$, many important features are different. For the latter
case, the method of moving spheres and the Pohozaev's type
arguments are very effective. However, these well known methods do
not seem to be useful for the former case. Instead we mainly use
the potential theory iteratively to obtain the sharp estimate.
Heuristically, the main difference between the two cases is that
when $\alpha>0 (\alpha\not\in \Bbb N)$, the linearized operator of
(\ref{mar28e2}) along a standard bubble is "invertible". This
invertibility forces the maximum points of blowup solutions to be
very close to $0$ (see Corollary \ref{nov10cor1}).  While for the
$\alpha=0$ case, this invertibility is lost and it is important to
obtain the sharp vanishing rate of the gradient of some
coefficient functions, see \cite{zhangcmp} for details.

We also note that the assumption $\alpha\not \in \mathbb N$ cannot
be removed in general. In fact even Theorem \ref{thm1} may not
hold for $\alpha \in \mathbb N$ (see \cite{BCLT}
\cite{tarantello2}). However the case $\alpha \in \mathbb N$ has
important applications (see \cite{tarantello}) and needs to be
better understood. A theorem similar to Theorem \ref{thm2} in this
case, even with more assumptions, should still be of much
interest.

The structure of this paper is as follows:  In section 2 we prove Theorem \ref{thm2}.
Our proof is based on Theorem \ref{thm1}.
Then in the appendix we include some detailed estimates.

\medskip

\noindent{\bf Acknowledgement} The author would like to thank the
referee for scrutinizing the paper and for his/her many insightful
remarks. He is also grateful to M. Lucia for stimulating
discussions.

\section{Proof of Theorem \ref{thm2}}
We shall always use $C$ to denote a constant depending on $\alpha,
C_0,C_1$ only, unless we specify otherwise.

\subsection{A uniqueness lemma}
\begin{lem}
\label{lem1}
 Let $\alpha\in \mathbb R^+ \setminus \mathbb N$,
$\phi$ be a $C^2$
 function that verifies
$$\left\{\begin{array}{ll}
\Delta \phi+|x|^{2\alpha}e^{U_{\alpha}}\phi=0\qquad
 \hbox{ in } \mathbb R^2,\\
\\
\phi(0)=0,\quad |\phi(x)|\le C(1+|x|)^{\tau}
 \quad x\in \mathbb R^2.
\end{array}
\right.
$$
where
$$
  U_{\alpha}(x)=\log \frac{8(\alpha+1)^2}{(1+|x|^{2\alpha+2})^2}
  \quad \hbox{ and } \quad
  \tau\in [0,1).
$$
Then $\phi \equiv 0$.
\end{lem}

\noindent{\bf Proof of Lemma \ref{lem1}:} Let $k\ge 1$ be an
integer. We define
$$\phi_k(r)=\frac 1{2\pi}\int_0^{2\pi}\phi(r\cos \theta, r\sin
 \theta)\cos (k\theta)d\theta.$$
Then $\phi_k$ satisfies
\begin{equation}
\left\{\begin{array}{ll} \phi_k''(r)+\frac
 1r\phi_k'(r)+(r^{2\alpha}e^{U_{\alpha}}-\frac{k^2}{r^2})\phi_k(r)=0, \quad
0<r<\infty,\\
\\
\lim_{r\to 0}\phi_k(r)=0,\quad |\phi_k(r)|\le C(1+r)^{\tau}, \quad
r>0.
\end{array}
\right. \label{nov9eq2}
\end{equation}
Putting $r=e^t$ and using Theorem 8.1 of \cite[Theorem 8.1,
p.92]{cod} one can see there are two fundamental solutions of
(\ref{nov9eq2}) that behave like those of
$$\tilde \phi''(r)+\frac{1}r\tilde \phi'(r)-\frac{k^2}{r^2}\tilde
 \phi(r)=0$$
as $r\to \infty$. (Note that the $\int_0^{\infty}|V'(t)|dt<\infty$
in the statement of Theorem 8.1 from \cite{cod} corresponds to
$\int_0^{\infty}r^{2\alpha+2}e^{U_{\alpha}}dr<\infty,$ which is
readily verified.)

Let us denote by $\phi_{k,1}$ and $\phi_{k,2}$ the two fundamental
solutions of (\ref{nov9eq2}) so that $\phi_{k,1}\sim r^k$,
$\phi_{k,2}\sim r^{-k}$ as $r\to \infty$. Since
$\phi_k=c_1\phi_{k,1}+c_2\phi_{k,2}$ and $|\phi_k|$ grows no
faster than $r^{\tau}$ as $r\to \infty$, we conclude that
 $c_1=0$. Hence $|\phi_k|\sim r^{-k}$ for $r$ large. On the other hand,
 when $r$ is close to $0$, the term
 $r^{2\alpha}e^{U_{\alpha}}$ is a perturbation
again, which means there are two fundamental solutions, say,
$\tilde \phi_{k,1}$ and $\tilde \phi_{k,2}$, comparable to $r^k$
and $r^{-k}$ as $r\to 0+$, respectively. Since $\phi_k(0)=0$, we
know that $|\phi_k(r)|\sim r^k$ as $r\to 0+$.

We claim that $\phi_k \equiv 0$ for all $k\ge 1$. To see this, let
$f_k(s)=\phi_k(s^{\frac 1{1+\alpha}})$, then we have
\begin{equation}
\label{nov16e1} f_k''(s)+\frac
 1sf_k'(s)+
 \left(
    \frac{8}{(1+s^2)^2}-\frac{k^2}{(1+\alpha)^2s^2}
 \right)f_k(s)=0,
\quad0<s<\infty.
\end{equation}
Using $\phi_k=O(r^{-k})$ at infinity and $\phi_k=O(r^{k})$ at
 $0$ we conclude that
\begin{equation}
\label{nov16e2} f_k(s)=O(s^{-\frac{k}{1+\alpha}})\quad
\mbox{at}\quad \infty, \qquad
f_k(s)=O(s^{\frac{k}{1+\alpha}})\quad \mbox{at}\quad 0.
\end{equation}
Let $\delta_1=\frac{k}{1+\alpha}$. Since $\alpha\not \in \Bbb N$
we see clearly that $\delta_1\neq 1$. Then by direct computation
we verify that the following two functions are two fundamental
solutions of (\ref{nov16e1}):
\begin{eqnarray}
f_{11}(s)=\frac{(\delta_1+1)s^{\delta_1}+(\delta_1-1)s^{\delta_1+2}}{1+s^2},
 \nonumber \\
f_{12}(s)=\frac{(\delta_1+1)s^{2-\delta_1}+(\delta_1-1)s^{-\delta_1}}{1+s^2}.
\label{jan1e1}
\end{eqnarray}
Consequently $f_k=c_1f_{11}+c_2f_{12}$ where $c_1$ and $c_2$ are
two constants. From (\ref{nov16e2}) we see that $c_1=0$ because
otherwise $|f_k|\sim s^{\delta_1}$ at $\infty$. Similarly by
observing the behavior of $f_k$ at $0$ we have $c_2=0$ because
otherwise $|f_k|\sim s^{-\delta_1}$ near $0$, which is a
contradiction to (\ref{nov16e2}) again. So $f_k\equiv 0$ for all
$k\ge 1$. This is
 equivalent to $\phi_k\equiv 0$ for
all $k\ge 1$. The same argument also shows that the projection of
$\phi$ along $\sin k\theta$ ($\forall k\ge 1$) is $0$. Finally,
from $\phi(0)=0$ we conclude $\phi(x)\equiv 0$.$\Box$

\subsection{Initial refinement}

In this subsection, we give the first improvement of Theorem
 \ref{thm1}. The main result
in this subsection is Proposition \ref{nov7p1}.

Let $\delta_i=e^{-\frac{u_i(0)}{2+2\alpha}}$ and we denote by
$z_i$ the maximum
 point of $u_i$. It is proved in
\cite{BCLT} that $\delta_i^{-1}z_i\to 0$. Note that we use
$u_i(0)$ to define $\delta_i$ because
$u_i(0)/u_i(z_i)=1+\circ(1)$. From the definition of $\psi_i$,
$u_i-\psi_i$ has no oscillation on $\partial B_1$ and it satisfies
$$
 \Delta (u_i-\psi_i)+|x|^{2\alpha}V_i(x)e^{(u_i-\psi_i)(x)}=0, \quad
 \hbox{ in } B_1$$
where $V_i(x)=H_i(x)e^{\psi_i(x)}$ is defined as before. Without
loss of generality we assume
$$V_i(0)=H_i(0)\to 8(\alpha+1)^2.$$
Let $x_i$ be the maximum point of $u_i-\psi_i$, then we still have
$\delta_i^{-1}x_i\to 0$ by the argument in
 \cite{BCLT}. Now we define $v_i$ as
$$
  v_i(y)=u_i(\delta_iy)-\psi_i(\delta_iy)-u_i(0),
  \quad y \in \Omega_i := B(0, \delta_i^{-1} ).
$$
We list some properties of $v_i$ implied by its definition:
$$\left\{\begin{array}{ll}
\Delta v_i(y)+|y|^{2\alpha}V_i(\delta_iy)e^{v_i(y)}=0,\quad  y\in
 \Omega_i:=B(0,\delta_i^{-1}),\\
\\
v_i(0)=0,\quad y_i:=\delta_i^{-1}x_i\to 0.\\
\\
v_i(y)\to (-2)\log (1+|y|^{2\alpha+2})\quad
\mbox{in}\quad C^2_{loc}(\mathbb R^2).\\
\\
v_i(y_1)=v_i(y_2),\quad \forall y_1,y_2\in \partial
\Omega_i:=\partial
 B(0,\delta_i^{-1}).
\end{array}
\right.
$$
Note that in the second equation above we use $y_i$ to denote the
maximum point of $v_i$. Let
$$U_i(y)=(-2)
  \log
 \left(1+\frac{V_i(0)}{8(\alpha+1)^2}|y|^{2\alpha+2} \right)$$
be a standard bubble which satisfies
\begin{equation}
\label{mar31e2} \Delta U_i+|y|^{2\alpha}V_i(0)e^{U_i(y)}=0, \quad
 \hbox{ in } \mathbb R^2.
\end{equation}
Then the conclusion of Theorem \ref{thm1} can be written as:
$$|v_i(y)-U_i(y)|\le C ,\quad |y|\le \delta_i^{-1}. $$

Let $w_i(y)=v_i(y)-U_i(y)$, then the following proposition is the
first improvement of Theorem \ref{thm1}:

\begin{prop}
\label{nov7p1}
 For any $\epsilon\in (0,\frac 12)$, there exists $C(\epsilon)>0$ such
 that for all
large $i$,
$$|w_i(y)|\le C\delta_i(1+|y|)^{\epsilon},\quad \hbox{ in } \Omega_i.$$
\end{prop}

\noindent{\bf Proof of Proposition \ref{nov7p1}:} We write the
equation for $w_i$ as
$$\left\{\begin{array}{ll}
\Delta
 w_i+r^{2\alpha}V_i(0)e^{\xi_i}w_i=O(\delta_i)(1+r)^{-3-2\alpha},\\ \\
w_i(0)=0,\quad |w_i(y)|\le C, \quad y\in \Omega_i,\quad
w_i|_{\partial \Omega_i}=\tilde a_i,
\end{array}
\right.
$$ where $\xi_i$ is obtained by the mean value theorem.
From Theorem \ref{thm1} we immediately have $\tilde a_i=O(1)$. Let
$$\tilde \Lambda_i=\max\frac{|w_i(y)|}{\delta_i(1+|y|)^{\epsilon}}\quad y\in
 \bar \Omega_i.$$
Our goal is to show $\tilde \Lambda_i=O(1)$. We prove this by a
contradiction. Suppose
 $\tilde \Lambda_i\to \infty$, then we use $\tilde y_i$ to denote a point where
 $\tilde \Lambda_i$ is assumed. Let
$$\bar w_i(y)=\frac{w_i(y)}{\tilde \Lambda_i\delta_i(1+|\tilde y_i|)^{\epsilon}}.$$

It readily follows from the definition of $\tilde \Lambda_i$ that
$\displaystyle{|\bar w_i(y)|\le
 \frac{(1+|y|)^{\epsilon}}{(1+|\tilde y_i|)^{\epsilon}}},$
which means $\bar w_i$ is uniformly bounded over any fixed compact
subset of $\Bbb R^2$. Therefore we conclude that, along a
subsequence, $\bar w_i$ converges in $C^2_{loc}(\mathbb R^2)$ to a
solution $w$ of
$$\left\{\begin{array}{ll}
\Delta w+r^{2\alpha}e^{U_{\alpha}}w=0,\quad \mathbb R^2,\\
\\
w(0)=0,\quad |w(y)|\le C(1+|y|)^{\epsilon}.
\end{array}
\right.
$$
If $\tilde y_i$ converges to $y_0\in \Bbb R^2$, we have
$|w(y_0)|=1$ by continuity. However this is impossible because an
application of Lemma \ref{lem1} shows that $w\equiv 0$. Therefore
the only case to consider is $\tilde y_i\to \infty$.

It follows from $|\bar w_i(\tilde y_i)|=1$ and the Green's
representation formula that
\begin{eqnarray}
\pm 1=\bar w_i(\tilde y_i)= \int_{\Omega_i}G(\tilde y_i,\eta)\bigg
\{|\eta
 |^{2\alpha}V_i(0)e^{\xi_i(\eta)}
\frac{w_i(\eta)}{\tilde \Lambda_i\delta_i(1+|\eta |)^{\epsilon}}
\frac{(1+|\eta |)^{\epsilon}}{(1+|\tilde y_i|)^{\epsilon}}\nonumber\\
+\frac{O(1)(1+|\eta |^{-3-2\alpha})}{\tilde \Lambda_i(1+|\tilde
y_i|)^{\epsilon}} \bigg \}d\eta -\int_{\partial
\Omega_i}\frac{\partial G}{\partial \nu} (\tilde y_i,\eta)
\frac{\tilde a_i}{\tilde \Lambda_i\delta_i(1+|\tilde
y_i|)^{\epsilon}}dS \label{jul17e1}
\end{eqnarray}
where $G$ is the Green's function over $\Omega_i$ with respect to
the Dirichlet boundary condition. Recall that the Green's function
over $\Omega_i$ is
$$G(y,\eta)=-\frac 1{2\pi}\log |y-\eta |
+\frac 1{2\pi}\log
 \left(
   \frac{|y|}{\delta_i^{-1}}|\frac{\delta_i^{-2}y}{|y|^2}-\eta |
 \right).$$

Since $\bar w_i(0)=0$, the Green's representation formula gives
\begin{eqnarray}
0=\int_{\Omega_i}G(0,\eta)\bigg
\{|\eta|^{2\alpha}V_i(0)e^{\xi_i(\eta)}
\frac{w_i(\eta)}{\tilde \Lambda_i\delta_i(1+|\eta |)^{\epsilon}}
\frac{(1+|\eta |)^{\epsilon}}{(1+|\tilde y_i|)^{\epsilon}}\nonumber\\
+\frac{O(1)(1+|\eta |)^{-3-2\alpha}}{\tilde \Lambda_i(1+|\tilde
y_i|)^{\epsilon}} \bigg \}d\eta -\int_{\partial
\Omega_i}\frac{\partial G}{\partial \nu} (0,\eta) \frac{\tilde
a_i}{\tilde \Lambda_i\delta_i(1+|\tilde y_i|)^{\epsilon}}dS.
\label{jul17e2}
\end{eqnarray}

To deal with the two boundary integral terms in (\ref{jul17e1})
and
 (\ref{jul17e2}),
we observe that
$$\int_{\partial \Omega_i}\bigg (\frac{\partial G}{\partial \nu}
(\tilde y_i,\eta) -\frac{\partial G}{\partial \nu}(0,\eta)\bigg )
\frac{\tilde a_i}{\tilde \Lambda_i\delta_i(1+|\tilde
y_i|)^{\epsilon}}dS =0.$$ This equality follows by using the well
known identity
$$\int_{\partial
\Omega_i}\partial_{\nu}G(\xi,\eta)dS=-1,\quad \forall \xi\in
\Omega_i.$$

From  (\ref{jul17e1}) and (\ref{jul17e2}) we have
\begin{equation}
\label{ccm1}
 1\le \int_{\Omega_i}|G(\tilde
y_i,\eta)-G(0,\eta)|\bigg (\frac{(1+|\eta
 |)^{-4-2\alpha+\epsilon}}
{(1+|\tilde y_i|)^{\epsilon}} +\circ(1)\frac{(1+|\eta
|)^{-3-2\alpha}}{(1+|\tilde y_i|)^{\epsilon}}\bigg ) d\eta.
\end{equation}
 Note that in the above we used
$$|\frac{w_i(\eta)}{\tilde \Lambda_i\delta_i(1+|\eta |)^{\epsilon}}|\le 1,
\quad e^{\xi_i(\eta)}\le C(1+|\eta |)^{-4-4\alpha},\quad \tilde
\Lambda_i\to \infty.$$

To get a contradiction to (\ref{ccm1}) we only need to show the
following:
\begin{equation} \label{ccm2}
\int_{\Omega_i}|G(\tilde y_i,\eta)-G(0,\eta)|\bigg (\frac{(1+|\eta
 |)^{-4-2\alpha+\epsilon}}
{(1+|\tilde y_i|)^{\epsilon}} +\circ(1)\frac{(1+|\eta
|)^{-3-2\alpha}}{(1+|\tilde y_i|)^{\epsilon}}\bigg )
d\eta=\circ(1).
\end{equation}
We consider two cases: If $|\tilde y_i|=\circ(1)\delta_i^{-1}$,
$G(\tilde y_i, \eta)$ can be written as
$$G(\tilde y_i,\eta)=-\frac 1{2\pi}\log |\tilde y_i-\eta |+\frac 1{2\pi}\log
 \delta_i^{-1}+\circ(1).$$
In this case it is enough to show
$$
\int_{\Omega_i} \bigg |\log \frac{|\tilde y_i-\eta |}{|\eta |}
\bigg | \bigg (\frac{(1+|\eta |)^{-4-2\alpha+\epsilon}}{(1+|\tilde
y_i|)^{\epsilon}} +\circ(1)\frac{(1+|\eta
|)^{-3-2\alpha}}{(1+|\tilde y_i|)^{\epsilon}}\bigg )
d\eta=\circ(1),$$ which follows from standard elementary
estimates.

Finally we consider the case $|\tilde y_i|\sim \delta_i^{-1}$. For
the Green's function we use
$$|G(\tilde y_i,\eta)-G(0,\eta)|\le C(\log (1+|\eta |)+\log \delta_i^{-1}),
\quad C \mbox{  universal }.$$ Then it is easy to obtain
(\ref{ccm2}) by elementary estimates. Proposition \ref{nov7p1} is
established. $\Box$

\medskip

From Proposition \ref{nov7p1} we obtain an estimate on $x_i$ more
precise than $|x_i|=\circ(\delta_i)$ (Recall that $x_i$ is the
place where the maximum of $u_i-\psi_i$ occurs, $y_i$ is the point
where the maximum of $v_i$ is attained.).
\begin{cor}
\label{nov10cor1}
\begin{equation}
\label{nov27e3}
|y_i|=|\delta_i^{-1}x_i|=O(\delta_i^{1/(2\alpha+1)}).
\end{equation}
\end{cor}
\noindent{\bf Proof of Corollary \ref{nov10cor1}:} Using
Proposition \ref{nov7p1} and standard elliptic estimates we have
$$|v_i(y)-U_i(y)|\le C\delta_i|y|,\quad |y|\le 10$$ and
$$|\nabla v_i(y)-\nabla U_i(y)|\le C\delta_i,\quad |y|\le 10.$$
Therefore $|\nabla U_i(y_i)|=O(\delta_i)$ because $\nabla
v_i(y_i)=0$. Then (\ref{nov27e3}) follows from the definition of
$U_i$. $\Box$

\begin{rem}
Recall that $z_i$ is the place where the maximum of $u_i$ is
attained. For $z_i$ we also have
$z_i=O(\delta_i^{\frac{2\alpha+2}{2\alpha+1}})$ by the same
argument.
\end{rem}

\subsection{Further refinement of the expansion}

$U_i$ is the first term in the expansion of $v_i$. To determine
the second term in the expansion we need a radial function $g_i$
that satisfies the following:
\begin{equation}
\label{nov10e7} \left\{\begin{array}{ll} g_i''(r)+\frac
1rg_i'(r)+(r^{2\alpha}V_i(0)e^{U_i(r)}-\frac
 1{r^2})g_i(r)
=-r^{2\alpha+1}e^{U_i(r)},\quad 0<r<\infty,\\
\\
\lim_{r\to 0+}g_i(r)=\lim_{r\to \infty}g_i(r)=0,\quad |g_i(r)|\le
 C\frac r{1+r^2}.
\end{array}
\right.
\end{equation}
By direct computation one checks that the following expression
verifies the above:
\begin{equation}
\label{mar27e1} g_i(r)=-\frac{2(1+\alpha)}{\alpha V_i(0)}
\frac{r}{1+\frac{V_i(0)}{8(1+\alpha)^2}r^{2\alpha+2}}.
\end{equation}

Let
\begin{equation}
\label{nov10e6}
\phi_i(y)=g_i(r)\delta_i\sum_{j=1}^2\partial_jV_i(0)\theta_j,\quad
 \theta_j=y_j/r.
\end{equation}
Then $\phi_i$ satisfies
\begin{equation}
\label{mar31e1} \Delta\phi_i+|y|^{2\alpha}V_i(0)e^{U_i(y)}\phi_i
=-\sum_{j=1}^2\delta_i\partial_jV_i(0)y_j|y|^{2\alpha}e^{U_i(y)},\quad
 \hbox{ in } \Omega_i.
\end{equation}

We claim that $\phi_i$ is the second term in the expansion of
$v_i$. This is verified in the following
\begin{prop}
\label{nov10p1} For any $\epsilon\in (0,\frac 12)$, there is $C>0$
depending only on
 $C_0,C_1,\alpha, \epsilon$
such that
\begin{eqnarray*}
&&|v_i(y)-(U_i+\phi_i)(y)|\le C\delta_i^2(1+|y|)^{\epsilon}, \\
&& |\nabla v_i(y)-\nabla (U_i+\phi_i)(y)|\le
 C\delta_i^2(1+|y|)^{\epsilon-1},\quad \hbox{ in } \Omega_i.
\end{eqnarray*}
\end{prop}

\noindent{\bf Proof of Proposition \ref{nov10p1}}: Let
$b_i=v_i-U_i-\phi_i$. Then $b_i$ satisfies $b_i(0)=0$ and
$$\Delta b_i+r^{2\alpha}V_i(\delta_iy)e^{\xi_i}b_i\\
=r^{2\alpha}(V_i(0)e^{U_i}-V_i(\delta_iy)e^{\xi_i})\phi_i
+O(\delta_i^2)r^{2\alpha+2}e^{U_i}.$$ By using Proposition
\ref{nov7p1} and (\ref{nov10e7}) in the estimate of the right hand
side of the above, we have
$$\Delta
 b_i+r^{2\alpha}V_i(\delta_iy)e^{\xi_i}b_i=O(\delta_i^2)(1+r)^{-2-2\alpha}.$$
Proposition \ref{nov10p1} follows from the same argument as in
Proposition \ref{nov7p1}. $\Box$.

\begin{rem}
The essential difference between the $\alpha\not\in\Bbb N$ case
and the $\alpha=0$ case lies in the fundamental solutions of the
equation in (\ref{nov10e7}). The two sets of fundamental solutions
behave very differently. As a result, for $\alpha\not\in \Bbb N$,
we have the presence of $\phi_i$ and the maximum point of $v_i$
very close to the singularity $0$ (see Corollary \ref{nov10cor1}).
On the other hand, for the $\alpha=0$ case, it is crucial to
obtain the vanishing rate of $\nabla V_i(0)$ from the Pohozaev
Identity (see \cite{zhangcmp}).
\end{rem}

Let $\theta_j := \frac{y_j}{r}$, $j=1,2$, we list the following
well known identities for convenience.
\begin{equation}
\label{aug23e2}
-\frac{d^2}{d\theta^2}(\theta_1\theta_2)=4\theta_1\theta_2,\quad
-\frac{d^2}{d\theta^2}(\theta_j^2-\frac 12)=4(\theta_j^2-\frac
 12),\quad j=1,2.
\end{equation}

Now we rewrite the equation for $v_i$ using Proposition
\ref{nov10p1}. Let $b_i=v_i-U_i-\phi_i$, then we have
\begin{equation}\label{ccm5}
\Delta (U_i+\phi_i+b_i)+r^{2\alpha}(V_i(0)+\delta_i\nabla
V_i(0)\cdot
 y+F_1+O(\delta_i^3r^3))
e^{U_i+\phi_i+b_i}=0
\end{equation}
 where
\begin{eqnarray}
F_1&=&\delta_i^2(\frac 12\partial_{11}V_i(0)y_1^2+\frac
 12\partial_{22}V_i(0)y_2^2+
\partial_{12}V_i(0)y_1y_2) \nonumber \\
&=&\delta_i^2\bigg (\frac
12\partial_{11}V_i(0)(y_1^2-\frac{r^2}2)+ \frac
 12\partial_{22}V_i(0)(y_2^2-\frac{r^2}2)+\partial_{12}V_i(0)y_1y_2\bigg )\nonumber \\
&&+\frac 14\delta_i^2r^2\Delta V_i(0)=F_{11}+F_{12}
\label{mar31e3}
\end{eqnarray}

\begin{eqnarray}
F_{11}&=&\delta_i^2r^2(\frac 12\partial_{11}V_i(0)(\theta_1^2-\frac 12)
+\frac 12 \partial_{22}V_i(0)((\theta_2^2-\frac 12)
+\partial_{12}V_i(0)\theta_1\theta_2) \nonumber\\
F_{12}&=&\frac 14\delta_i^2r^2\Delta V_i(0)
\label{mar31e4}
\end{eqnarray}

By using Proposition \ref{nov10p1} we have
$|b_i|=O(\delta_i^2)r^{\epsilon}$ and
$b_i+\phi_i=O(\delta_i(1+r)^{-\frac 12})$. Since
 $e^{U_i}=O((1+r)^{-4-4\alpha})$ we conclude that
\begin{eqnarray*}
e^{U_i+\phi_i+b_i}&=&e^{U_i}(1+\phi_i+b_i+\frac
 12(\phi_i+b_i)^2+O(\delta_i^3)r^{-\frac 32})\\
&=&e^{U_i}+e^{U_i}\phi_i+e^{U_i}b_i+\frac
 12e^{U_i}\phi_i^2+O(\delta_i^3)(1+r)^{-5-4\alpha+2\epsilon}.
\end{eqnarray*}

With this expression (\ref{ccm5}) is reduced to
\begin{eqnarray*}
&&\Delta(U_i+\phi_i+b_i)+r^{2\alpha}\bigg (V_i(0)+\delta_i\nabla
 V_i(0)\cdot y
+F_1+O(\delta_i^3r^3)\bigg )\\
&&\cdot \bigg (e^{U_i}+e^{U_i}\phi_i+e^{U_i}b_i+\frac
12e^{U_i}\phi_i^2 +O(\delta_i^3(1+r)^{-5-4\alpha+2\epsilon})\bigg )=0.
\end{eqnarray*}

The last term on the left hand side of the above is
\begin{eqnarray*}
&&r^{2\alpha}\bigg
 (V_i(0)e^{U_i}+V_i(0)e^{U_i}\phi_i+V_i(0)e^{U_i}b_i+\frac 12
V_i(0)e^{U_i}\phi_i^2+\delta_i\nabla V_i(0)\cdot ye^{U_i}\\
&&\quad +\delta_i\nabla V_i(0)\cdot ye^{U_i}\phi_i+F_1e^{U_i}
+O(\delta_i^3)(1+r)^{-1-4\alpha}\bigg ).
\end{eqnarray*}

By combining (\ref{mar31e1}), (\ref{mar31e2}), (\ref{mar31e3}) we
can write the equation for $b_i$ as
\begin{eqnarray}
&&\Delta
 b_i+r^{2\alpha}V_i(0)e^{U_i}b_i+\frac{V_i(0)}2r^{2\alpha}e^{U_i}\phi_i^2
+\delta_i\nabla V_i(0)\cdot ye^{U_i}\phi_ir^{2\alpha}\nonumber \\
&&+\frac{r^{2+2\alpha}}4\Delta
 V_i(0)\delta_i^2e^{U_i}+r^{2\alpha}F_{11}e^{U_i}
+O(\delta_i^3)(1+r)^{-1-2\alpha}=0. \label{mar19e1}
\end{eqnarray}

Without loss of generality we assume $\nabla V_i(0)=|\nabla
 V_i(0)|e_1$. Then the sum of the third and fourth term on
 the left hand side of
 (\ref{mar19e1}) is
\begin{eqnarray*}
&&\frac{V_i(0)}2r^{2\alpha}e^{U_i}\phi_i^2+\delta_ir^{2\alpha}\nabla
 V_i(0)\cdot ye^{U_i}\phi_i\\
&=&\frac{V_i(0)}2r^{2\alpha}e^{U_i}g^2_i(r)\delta_i^2|\nabla
 V_i(0)|^2\theta_1^2
+\delta_i^2g_i(r)|\nabla V_i(0)|^2r^{1+2\alpha}\theta_1^2e^{U_i}\\
&=&\delta_i^2r^{2\alpha}e^{U_i}|\nabla
 V_i(0)|^2\theta_1^2(\frac{V_i(0)}2g_i^2(r)
+g_i(r)r)\\
&=&\delta_i^2r^{2\alpha}e^{U_i}|\nabla V_i(0)|^2(\theta_1^2-\frac
12)
(\frac{V_i(0)}2g_i^2(r)+g_i(r)r)\\
&&\quad +\frac 12 \delta_i^2r^{2\alpha}e^{U_i}|\nabla V_i(0)|^2
(\frac{V_i(0)}2g_i^2(r)+g_i(r)r)\\
&=&C_{11}+C_{12}.
\end{eqnarray*}
For $C_{11}$ and $r^{2\alpha}F_{11}e^{U_i}$ we can find $c_i$ that
satisfies
\begin{equation}
\label{mar27e5} \left\{\begin{array}{ll} \Delta
 c_i+r^{2\alpha}V_i(0)e^{U_i}c_i+C_{11}+r^{2\alpha}F_{11}e^{U_i}=0,
\quad 0<r<\delta_i^{-1}\\ \\
|c_i(x)|\le C\delta_i^2\frac{r^2}{(1+r)^3},\quad
0<r<\delta_i^{-1}.
\end{array}
\right.
\end{equation}
The existence of $c_i$ and its estimate are established in the
appendix. Let $d_i=b_i-c_i$. Then the equation for $d_i$ is
\begin{equation}
\label{mar19e2} \Delta
 d_i+r^{2\alpha}V_i(0)e^{U_i}d_i+E_i+O(\delta_i^3(1+r)^{-1-2\alpha})=0.
\end{equation}
where
$$E_i=\frac{r^{2+2\alpha}}4\Delta V_i(0)\delta_i^2e^{U_i}
+\frac 12\delta_i^2r^{2\alpha}e^{U_i}|\nabla
 V_i(0)|^2(\frac{V_i(0)}2g_i^2(r)+g_i(r)r).$$

It follows from the definition of $d_i$, Proposition \ref{nov10p1}
and (\ref{mar27e5}), that
\begin{equation}
\label{mar27e3} d_i(0)=0,\quad |d_i(y)|\le
C\delta_i^2(1+r)^{\epsilon},\quad r\le
 \delta_i^{-1}.
\end{equation}
Also, it is implied by (\ref{mar19e2}) and (\ref{mar27e3}) that
\begin{equation}
\label{mar31e5} \int_{\Omega_i}\Delta d_i=O(\delta_i^2).
\end{equation}
In the following proposition we evaluate the value of $d_i$ on
$\partial \Omega_i$:
\begin{prop}
\label{mar27p1}
$$d_i=(\Lambda_1\Delta V_i(0)+\Lambda_2|\nabla V_i(0)|^2)\delta_i^2\log
 \delta_i^{-1}
+O(\delta_i^2)\quad \mbox{on}\quad \partial \Omega_i.$$
\end{prop}

\noindent{\bf Proof of Proposition \ref{mar27p1}:} Let
$$
f_i(y)=\frac{1-a_ir^{2\alpha+2}}{1+a_ir^{2\alpha+2}},\quad
 a_i=\frac{V_i(0)}{8(1+\alpha)^2}.
$$
Direct computation shows
\begin{equation}
\label{mar31e6} \Delta
f_i(y)+r^{2\alpha}V_i(0)e^{U_i}f_i(y)=0,
 \quad \hbox{ in } \mathbb R^2.
\end{equation}
Also it is straightforward to verify that
\begin{equation}
\label{mar20e1} f_i(y)=-1+O(\delta_i^{2\alpha+2}),\quad
\frac{\partial f_i}{\partial r}=O(\delta_i^{3+2\alpha})
 \quad \mbox{ on } \partial \Omega_i.
\end{equation}

A direct consequence of (\ref{mar19e2}), (\ref{mar31e6}) and the
Green's formula is
\begin{equation}
\label{mar20e2} \int_{\partial \Omega_i}(\frac{\partial
f_i}{\partial \nu}d_i -\frac{\partial d_i}{\partial
\nu}f_i)=\int_{\Omega_i}E_if_i +O(\delta_i^{2+2\alpha}).
\end{equation}

From (\ref{mar20e1}) we have
\begin{equation}
\label{mar20e3} \int_{\partial \Omega_i}\frac{\partial
f_i}{\partial \nu}d_i =O(\delta_i^{2+2\alpha}).
\end{equation}
Since $f_i$ is a radial function, we have, by (\ref{mar31e5}) and
 (\ref{mar20e1}),
\begin{equation}
\label{mar20e4} \int_{\partial \Omega_i}\frac{\partial
d_i}{\partial \nu}f_i =(\int_{\Omega_i}\Delta d_i)f_i|_{\partial
\Omega_i} =-\int_{\Omega_i}\Delta d_i+O(\delta_i^{2+2\alpha}).
\end{equation}

So we conclude from (\ref{mar20e2}), (\ref{mar20e3}) and
 (\ref{mar20e4}) that
\begin{equation}
\label{mar20e5} \int_{\Omega_i}\Delta
 d_i=\int_{\Omega_i}E_if_i+O(\delta_i^{2+2\alpha}).
\end{equation}
On the other hand, from $d_i(0)=0$ and (\ref{mar19e2}) we have
$$
0=\int_{\Omega_i}G(0,\eta)(|\eta
 |^{2\alpha}V_i(0)e^{U_i(\eta)}d_i(\eta)
+E_i(\eta))+d_i|_{\partial \Omega_i}+O(\delta_i^2).$$ Since
$$G(0,\eta)=-\frac{1}{2\pi}\log |\eta |+\frac{1}{2\pi}\log
 \delta_i^{-1}$$
and
$$\int_{\Omega_i}\log |\eta |(|\eta
 |^{2\alpha}V_i(0)e^{U_i(\eta)}d_i(\eta)+E_i(\eta))
=O(\delta_i^2),$$ we obtain from elementary estimates
\begin{eqnarray*}
 d_i|_{\partial \Omega_i}&=&\frac 1{2\pi}\log
 (\delta_i^{-1})\int_{\Omega_i}E_if_i+O(\delta_i^2)\\
&=&(\log \delta_i^{-1})\delta_i^2(\Lambda_1\Delta
V_i(0)+\Lambda_2|\nabla
 V_i(0)|^2)+
O(\delta_i^{2}).
\end{eqnarray*}
Proposition \ref{mar27p1} is established. $\Box$

\bigskip

We now finish the proof of Theorem \ref{thm2} by a standard
application of the maximum principle. Let
$$M_i(y)=(\Lambda_1\Delta V_i(0)+\Lambda_2|\nabla
 V_i(0)|^2)\delta_i^2\log |y|
+\delta_i^2M(1-r^{-\alpha})$$ where $M$ is a large number to be
 determined. Thanks to (\ref{mar19e2}),(\ref{mar27e3}) and the
 estimate of $E_i$, the equation for $d_i$ can be written as
$$\Delta d_i=O(\delta_i^2)(1+r)^{-2-2\alpha},
 \quad \hbox{ in } \Omega_i.$$
By choosing $M$ large enough we have $M_i(y)>d_i(y)$ on $\partial
 B_2$ and
$\partial \Omega_i$. On the other hand
$$\Delta M_i(y)=-M\alpha^2\delta_i^2r^{-2-\alpha}<\Delta d_i(y)\quad
 2<r<\delta_i^{-1}.$$
Therefore by the maximum principle
$$d_i(y)\le M_i(y)\quad 2<r<\delta_i^{-1}.$$ Similarly we can also show
$$d_i(y)\ge (\Lambda_1\Delta V_i(0)+\Lambda_2|\nabla
 V_i(0)|^2)\delta_i^2\log |y|
-\delta_i^2M(1-r^{-\alpha})\quad 2<r<\delta_i^{-1}$$ for $M$
large. Theorem \ref{thm2} is established. $\Box$

\section{Appendix}

In this section we prove the existence of $c_i$ that satisfies
(\ref{mar27e5}).
 Recall that
$$F_{11}=\delta_i^2r^2(\frac 12\partial_{11}V_i(0)(\theta_1^2-\frac 12)
+\frac 12\partial_{22}V_i(0)(\theta_2^2-\frac
 12)+\partial_{12}V_i(0)\theta_1\theta_2)$$
$$C_{11}=\delta_i^2r^{2\alpha}e^{U_i}|\nabla V_i(0)|^2(\theta_1^2-\frac
 12)
(\frac{V_i(0)}2g_i^2(r)+g_i(r)r).$$ The three terms in
$r^{2\alpha}F_{11}e^{U_i}$ and $C_{11}$ can all be
 written in the form
$\delta_i^2f(\theta)Q_i(r)$ where $f(\theta)$ is one of the
spherical harmonics ( $\theta_1^2-\frac 12$, $\theta_2^2-\frac
12$,
 $\theta_1\theta_2$ ) and
 $Q_i(r)$ is a radial
function that satisfies
$$|Q_i(r)|\le C\frac{r^{2+2\alpha}}{(1+a_ir^{2+2\alpha})^2}.$$
If we can find a function $h_1(r)$ that solves
\begin{equation}
\label{mar27e6} \left\{\begin{array}{ll} h_1''(r)+\frac
1rh_1'(r)+(r^{2\alpha}V_i(0)e^{U_i}-\frac 4{r^2})h_1(r)
=-Q_i(r),\quad 0<r<\delta_i^{-1}\\ \\
|h_1(r)|\le C\frac{r^2}{(1+r)^3},\quad 0<r<\delta_i^{-1},
\end{array}
\right.
\end{equation}
then by (\ref{aug23e2}) $h_1(r)f(\theta)$ solves
$$(\Delta+r^{2\alpha}V_i(0)e^{U_i})(f(\theta)h_1(r))+f(\theta)Q_i(r)=0,\quad
\Omega_i.$$ Consequently $c_i$ is the sum of four such functions.
So all we need to establish is (\ref{mar27e6}). Let
$$r=(\sqrt{\frac{8(\alpha+1)^2}{V_i(0)}}s)^{\frac 1{1+\alpha}}\quad
 \mbox{ and }
f_1(s)=h_1\bigg ((\sqrt{\frac{8(\alpha+1)^2}{V_i(0)}}s)^{\frac
 1{1+\alpha}}\bigg ).$$
Then the equation for $f_1(s)$ is
$$f_1''(s)+\frac 1sf_1'(s)+
(\frac 8{(1+s^2)^2}-\frac{4\delta^2}{s^2})f_1(s)=l_1(s),\quad
 0<s<\infty,$$
where $\delta=\frac 1{1+\alpha}$,
$$|l_1(s)|\le C\frac{s^{2\delta}}{(1+s^2)^2}.$$
Since $\alpha$ is not an integer, $2\delta-1\neq 0$, two
fundamental solutions of the homogeneous equation ($f_{21}$ and
$f_{22}$) can be found in explicit form:
\begin{eqnarray*}
f_{21}(s)=\frac{(2\delta+1)s^{2\delta}+(2\delta-1)s^{2\delta+2}}{1+s^2}\\
f_{22}(s)=\frac{(2\delta+1)s^{2-2\delta}+(2\delta-1)s^{-2\delta}}{1+s^2}
\end{eqnarray*}
Let
$$w_1(s):=f_{21}(s)f_{22}'(s)-f'_{21}(s)f_{22}(s)=4\delta(1-4\delta^2)s^{-1}$$
and
$$
f_1(s)=-\int_s^{\infty}\frac{f_{22}(\tau)l_1(\tau)}{w_1(\tau)}d\tau
 f_{21}(s)
+\int_0^s\frac{f_{21}(\tau)l_1(\tau)}{w_1(\tau)}d\tau f_{22}(s),
$$
then it is straightforward to verify
$$|f_1(s)|\le C\frac{s^{2\delta}}{(1+s)^3},\quad 0<s<\infty.$$
(\ref{mar27e6}) is established. $\Box$

\end{document}